\documentclass[pdflatex,sn-mathphys-num]{sn-jnl}

\usepackage{graphicx}%
\usepackage{multirow}%
\usepackage{amsmath,amssymb,amsfonts}%
\usepackage{amsthm}%
\usepackage{mathrsfs}%
\usepackage{algorithm}
\usepackage[title]{appendix}%
\usepackage{xcolor}%
\usepackage{textcomp}%
\usepackage{manyfoot}%
\usepackage{booktabs}%
\usepackage{algorithm}%
\usepackage{algorithmicx}%
\usepackage{algpseudocode}%
\usepackage{listings}%

\theoremstyle{thmstyleone}%
\newtheorem{theorem}{Theorem}
\newtheorem{proposition}[theorem]{Proposition}%
\newtheorem{assumption}{Assumption}
\newtheorem{lemma}{Lemma}
\theoremstyle{thmstyletwo}%
\newtheorem{remark}{Remark}%

\theoremstyle{thmstylethree}%
\numberwithin{equation}{section}
\numberwithin{theorem}{section}
\begin{document}


\title[ M-LBFGS with displacement aggregation and its application]{Modified limited memory BFGS with displacement aggregation and its application to the largest eigenvalue problem}


\author*[1]{\fnm{Manish Kumar} \sur{Sahu}}\email{manishkumarsahu132@gmail.com}

\author[2]{\fnm{Suvendu Ranjan} \sur{Pattanaik}}\email{suvendu.pattanaik@gmail.com}

\affil[]{\orgname{National Institute of Technology Rourkela, India}} 


\abstract{  
We present a modified limited memory BFGS method with displacement aggregation (AggMBFGS) for solving nonconvex optimization problems. AggMBFGS refines curvature pair updates by removing linearly dependent variable variations, ensuring that the inverse Hessian approximation retains essential curvature properties. As a result, its per iteration complexity and storage requirement is $\mathcal{O}(\tau d)$ where  $\tau \leq d$ represents the memory size and $d$ is the problem dimension. We establish the global convergence of both M-LBFGS and AggMBFGS under a backtracking modified Armijo line search (MALS) and prove the local superlinear convergence of AggMBFGS, demonstrating its theoretical advantages over M-LBFGS with the classical Armijo line search~\cite{Shi2016ALM}. Numerical experiments on CUTEst test problems~\cite{gould2015cutest} confirm that AggMBFGS outperforms M-LBFGS in reducing the number of iterations and function evaluations. Additionally, we apply AggMBFGS to compute the largest eigenvalue of high-dimensional real symmetric positive definite matrices, achieving lower relative errors than M-LBFGS~\cite{Shi2016ALM} while maintaining computational efficiency. These results suggest that AggMBFGS is a promising alternative for large-scale nonconvex optimization and eigenvalue computation.  
}

\keywords{Displacement aggregation (DA), Modified  Broyden Fletcher Goldfarb Shanno (MBFGS), Modified Limited Memory Broyden Fletcher Goldfarb Shanno (M-LBFGS),  Modified Limited Memory Broyden Fletcher Goldfarb Shanno with displacement aggregation(AggMBFGS), Modified Armijo Line Search (MALS)}



\maketitle

\section{Introduction}\label{sec1}

In this paper, we consider the unconstrained optimization problem:
\begin{equation} \label{eq:main_problem}
	\min_{x \in \mathbb{R}^{d}} f(x),
\end{equation}
where \( f: \mathbb{R}^{d} \to \mathbb{R} \) is a twice continuously differentiable function. One of the most effective approaches to solving \eqref{eq:main_problem} is the quasi-Newton method, which iteratively approximates the inverse Hessian matrix using gradient information. The quasi-Newton method often outperforms both gradient descent and Newton's method in terms of computational efficiency and storage requirements \cite{nocedal1999numerical}.

Newton's method exhibits quadratic convergence \cite{nocedal1999numerical}; however, explicitly computing and storing the inverse Hessian matrix is computationally expensive, particularly for large-scale problems. To address this issue, quasi-Newton methods approximate the inverse Hessian using various update strategies such as Symmetric Rank-1 (SR1), Davidson Fletcher Powell (DFP), and Broyden Fletcher Goldfarb Shanno (BFGS) \cite{nocedal1999numerical}. Among these, the BFGS update is known for its superior empirical performance. While quasi-Newton methods achieve superlinear convergence with reduced computational overhead, their full-memory versions remain impractical for high-dimensional problems.

To overcome this limitation, Liu and Nocedal \cite{Liu1989OnTL} introduced the limited-memory BFGS (L-BFGS) method, which maintains a restricted number of curvature pairs, controlled by a memory parameter \( \tau \). Unlike BFGS, L-BFGS avoids explicit inverse Hessian computation by leveraging recursive two-loop updates. Although L-BFGS is computationally efficient, it typically exhibits only linear convergence.

Despite its widespread success, standard BFGS methods may fail to converge for nonconvex functions \cite{Dai2002ConvergencePO, Mascarenhas2004TheBM}. Consequently, several modifications have been proposed to enhance stability and convergence in nonconvex settings \cite{Li2001AMB,babaie2011modified}. Notably, Li and Fukushima \cite{Li2001AMB} developed a modified BFGS (MBFGS) update that ensures global convergence and achieves superlinear convergence under suitable conditions. Xiao et al. \cite{Xiao2013GlobalCO} extended this approach, proving the global convergence of M-LBFGS with Wolfe line search for nonconvex functions. Similarly, Shi et al. \cite{Shi2016ALM} established global convergence guarantees for M-LBFGS under an Armijo line search framework and applied M-LBFGS to the largest eigenvalue problem. However, extending the theoretical convergence properties of full-memory MBFGS to its limited-memory variant remains challenging due to reduced curvature information and storage constraints. 

In 2022, Berahas and Curtis \cite{Berahas2019LimitedmemoryBW} introduced the AggBFGS method, which employs displacement aggregation to eliminate linearly dependent variable variations while refining curvature information. They demonstrated that AggBFGS retains the convergence properties of full-memory BFGS while requiring fewer iterations and function evaluations than standard L-BFGS. Despite these advancements, achieving superlinear convergence in limited-memory quasi-Newton methods for nonconvex optimization remains an open problem.

Given these limitations, we formulate two key research questions:
\begin{enumerate}
	\item Can we design a limited memory quasi-Newton method that ensures superlinear convergence while maintaining low computational complexity and memory efficiency in nonconvex settings?
	\item What are the key application areas where such an improved limited-memory quasi-Newton method would provide significant benefits, particularly in high-dimensional optimization problems?
\end{enumerate}
\subsection*{Contributions}  
\begin{enumerate}  
	\item We introduce the AggMBFGS method, a modified limited-memory BFGS algorithm incorporating displacement aggregation. This approach refines the updates of curvature pairs by removing linearly dependent variable variations, ensuring a more accurate inverse Hessian approximation while maintaining a memory complexity of $\mathcal{O}(\tau d)$.
	
	\item We establish the local superlinear rate of convergence and global convergence of AggMBFGS under a backtracking modified Armijo line search. Theoretical analysis confirms that AggMBFGS provides improved convergence guarantees over M-LBFGS, particularly for nonconvex optimization problems. Numerical experiments on CUTEst test problems \cite{gould2015cutest} validate its effectiveness.
	
	\item We apply AggMBFGS with a modified Armijo line search to compute the largest eigenvalue of high-dimensional real symmetric positive definite matrices. The numerical results show that AggMBFGS achieves lower relative errors than M-LBFGS \cite{Shi2016ALM} while maintaining computational efficiency, making it a competitive alternative for large-scale eigenvalue problems.
	
\end{enumerate}


\par The remainder of this paper is structured as follows. Section \ref{293} introduces preliminaries on the relationship between the minimization of a function \( f(x) \) and the eigenvalues of a real positive definite matrix \( A \). Additionally, it provides a brief overview of MBFGS, M-LBFGS, and the modified Armijo line search (MALS). In Section \ref{789}, we propose AggMBFGS and establish its local superlinear rate of convergence and global convergence with the modified Armijo line search. We further demonstrate that AggMBFGS inherits the convergence properties of full-memory MBFGS and outperforms M-LBFGS with the modified Armijo line search. 
Section~\ref{3456} applies the proposed AggMBFGS method with a modified Armijo line search to compute the largest eigenvalues of sparse matrices from the University of Florida collection~\cite{davis2011university}. It then analyzes the error in estimating the largest eigenvalue compared to the \texttt{eigs} command in MATLAB. Finally, Section~\ref{5876} presents our conclusions.
\section{Preliminaries}\label{293}
A point \( x^* \) is called a critical point of a differentiable function \( f(x) \) if its gradient vanishes, i.e., \( g(x) = \nabla f(x) = 0 \). A critical point is classified as a local minimum, local maximum, or saddle point depending on the definiteness of the Hessian matrix at that point: it is positive definite for a local minimum, negative definite for a local maximum, and indefinite for a saddle point \cite{Auchmuty1989UnconstrainedVP}.  

\par In 1989, Auchmuty \cite{Auchmuty1989UnconstrainedVP} demonstrated that the largest eigenvalues of a real symmetric positive definite matrix \( A \)  can be computed by solving the following optimization problem:   
\par \begin{equation}\label{420}
	\min_{ x \in \mathbb{R}^{d}} f(x)=\min_{ x \in \mathbb{R}^{d}}  \frac{\|x \|^4}{4}- \frac { x^{T}Ax}{2}.
\end{equation}

\begin{theorem}
	( \cite{Auchmuty1989UnconstrainedVP},Theorem 12)  Let $f(x)$ be defined as (\ref{420}), a twice continuously differentiable, non-convex function, and  A is a real symmetric positive definite matrix. Then
	\begin{enumerate}
		\item  $f(x)$ is coercive and $\min_{ x \in \mathbb{R}^{d}} f(x)=-\frac{\lambda_1^2}{4} $ and $\lambda_{1}$ is the largest eigenvalue of A.
		\item  $x^*=\sqrt{\lambda_{j}}e_j$  is the non-zero critical point of f where $\lambda_{j}$ is a positive eigenvalue of A and $e_j$ is a normalized eigenvector corresponds to $\lambda_{j}$. Moreover, if $\lambda_{1} \neq \lambda_{j}$,$\sqrt{\lambda_{j}}e_j$ is a saddle point.
	\end{enumerate} 
\end{theorem}
Several optimization methods have been employed to compute the largest eigenvalue of real positive definite matrices, including  Barzilai-Borwein  method~\cite{Gao2014BarzilaiBorweinlikeMF}, gradient descent method, Newton's method, and quasi-Newton methods~\cite{Mongeau2004ComputingEO}. Among these, the quasi-Newton method is the most effective for computing the largest eigenvalue of real positive definite matrices. Zhanwen et al. \cite{Shi2016ALM} demonstrated that M-LBFGS with Armijo line search can solve the largest eigenvalue problem with moderate accuracy. However, the step length obtained via the Armijo line search does not always ensure a significant decrease in the objective function \( f(x) \) \cite{Armijo1966MinimizationOF}. 
\begin{remark}
	Minimizing the function $
	f(x) = \frac{\|x \|^4}{4} - \frac{x^{T} A x}{2}
	$
	using an optimization algorithm yields the largest eigenvalue of the matrix \( A \) defined in \( f(x) \). The BFGS method is one of the most effective quasi-Newton method for solving unconstrained convex optimization problems. However, since \( f(x) \) is a twice continuously differentiable (\( C^2 \)) but non-convex function, BFGS may not be a suitable choice. In particular, BFGS does not guarantee convergence for certain nonconvex functions \cite{Mascarenhas2004TheBM}.
\end{remark}
\subsection{Background on MBFGS and M-LBFGS}\label{435}
We address problem~\eqref{eq:main_problem} using a modified quasi-Newton algorithm that approximates the inverse Hessian to avoid the $\mathcal{O}(d^3)$ complexity of direct computation in Newton's method. The update rule is given by  
\begin{equation}
	x_{t+1} = x_t - \alpha_t W_t \nabla f(x_t),
\end{equation}
where $\alpha_t$ is the step length, $x_t$ is the decision variable, $\nabla f(x_t)$ is the gradient, and $W_t$ is the inverse Hessian approximation at iteration $t$. 

In MBFGS \cite{Li2001AMB}, given $G_t \approx \nabla^2 f(x_t)$, we define the variable variation $s_t := x_{t+1} - x_t$, $y_t := \nabla f(x_{t+1}) - \nabla f(x_t)$, the gradient $g_t = \nabla f(x_t)$, and 
\[
\bar{y}_t = y_t + r_t \|g_t\| s_t, \quad \text{where} \quad r_t = 1 + \max\left(0, -\frac{y_t^T s_t}{s_t^T s_t}\right).
\]
The updated Hessian approximation $G_{t+1} \approx \nabla^2 f(x_{t+1})$ is chosen to be close to $G_t$ while satisfying the secant condition $G_{t+1} s_t = \bar{y}_t$. This yields the MBFGS update:
\begin{equation}\label{9987}
	G_{t+1} = G_t + \frac{\bar{y}_t \bar{y}_t^T}{\bar{y}_t^T s_t} - \frac{G_t s_t s_t^T G_t}{s_t^T G_t s_t} =: BFGS(G_t, s_t, y_t).
\end{equation}  
Setting $H_t = G_t^{-1}$, the rank-two MBFGS update allows efficient computation of $W_{t+1}$ via the Sherman-Morrison formula:
\begin{equation}\label{1111}
	W_{t+1} = \left(I - \frac{\bar{y}_t s_t^T}{s_t^T \bar{y}_t}\right)^T W_t \left(I - \frac{\bar{y}_t s_t^T}{s_t^T \bar{y}_t}\right) + \frac{s_t s_t^T}{s_t^T \bar{y}_t}.
\end{equation}
The memory and computational cost of the MBFGS update is $\mathcal{O}(d^2)$.  

Alternatively, MBFGS can be implemented by storing curvature pairs $(s_i, \bar{y}_i)$ for all $1 \leq i \leq t$ and computing $W_t$ from scratch at each iteration. The modified limited-memory BFGS (M-LBFGS) retains only the most recent $\tau$ curvature pairs $\{(s_i, \bar{y}_i)\}_{t-\tau+1}^{t}$. The descent direction $d_t = -W_t \nabla f(x_t)$ is then efficiently computed using the two-loop recursion ( see Algorithm 7.4 in \cite{nocedal1999numerical}), with a per iteration complexity of $\mathcal{O}(\tau d)$. For $\tau < d$, M-LBFGS has lower computational and memory costs than MBFGS.

\begin{algorithm}
	\caption{M-LBFGS with backtracking modified Armijo line search}\label{506}
	\begin{algorithmic}[1]
		\Require an initial guess $x_0 \in \mathbb{R}^n$ and  $W_0    \succ 0$ be the initial inverse Hessian approximation and constants $\sigma \in (0,1)$,  $\mu \in [0,\infty)$, $L_0>0$, $m>0$,  $p\in(0,1)$ and $\epsilon>0$
		\State Let $t=0$.
		\State Find $g_t=\nabla f(x_t)$.
		\If{$\|g_t\| < \epsilon$}
		\State stop.
		\Else
		\State Let $\tau := \min(t+1,\tau)$.
		\State Find $s_{t}=x_{t+1}-x_{t},\quad y_{t}=g_{t+1} -g_{t},\quad r_{t}=1+ \max[0,-\frac {y_{t}^Ts_t}{s_t^T s_t}]$.
		\State Compute $\bar{y}_{t}=y_{t}+r_{t}\|g_{t}\|s_{t} $.
		\State Update $x_{t+1}=x_{t}+\alpha_{t}d_{t}$ where $\alpha_t$ is computed using Algorithm \ref{alg:2} and $d_t$ is computed using Algorithm \ref{alg3}.
		\State Then  $t=t+1$ and go to step-2.
		\EndIf
	\end{algorithmic}
\end{algorithm}
\subsection{Modified Armijo line search (MALS)}
The classical Armijo line search \cite{Armijo1966MinimizationOF} is defined as 
\begin{equation}
	\label{7890}
	f(x_t + \alpha_{t}d_t) \leq f(x_t) + \sigma \alpha_{t} g_t^{T}d_t, 
\end{equation}
where $g_t=\nabla f(x_t)$, $\sigma \in (0,1)$, $\alpha_t = p^{j_t}$, $p \in (0,1)$, and $j_t$ is the smallest nonnegative integer such that it satisfies (\ref{7890}). However, in the classical Armijo line search, the reduction in the objective function may be small. To address this limitation, Zhong Whan proposed a modified Armijo line search \cite{Wan2012NewCB}, which is defined as  
\begin{equation}
	\label{7777}
	f(x_t + \alpha_{t}d_t) \leq f(x_t) + \sigma \alpha_{t} \left( g_t^{T}d_t - \alpha_{t} \mu L_t \| d_t \|^2 \right),
\end{equation}
where $\beta_t = -g_{t}^T d_t / (L_{t} \|d_t \|^2)$, and $\alpha_{t}$ is the largest component of the sequence $\{\beta_t, \beta_t p, \beta_t p^{2},\dots \}$ that satisfies (\ref{7777}). Here, $p \in (0,1)$, $\mu \in [0,\infty)$, $\sigma \in (0,1)$, and $L_t$ is an approximation to the Lipschitz constant $L$. A cautious BFGS method combined with the modified Armijo line search has been shown to perform better than one using the classical Armijo line search \cite{Wan2012NewCB}.

\begin{algorithm}
	\caption{Backtracking Modified Armijo Line Search}
	\label{alg:2}
	\begin{algorithmic}[1]
		\Require Initial point \( x_0 \), search direction \( d_t \), initial step size \( \alpha_0 > 0 \), parameters \( \sigma \in (0,1) \), \( \mu \geq 0 \), \( p \in (0,1) \), and initial estimate \( L_0 > 0 \).  
		\If{$t > 1$}
		\State Update \( L_t = \frac{s_{t-1}^T \bar{y}_{t-1}}{\|s_{t-1}\|^2} \)
		\EndIf
		\State Compute \( \beta_t = -\frac{g_t^T d_t}{L_t \| d_t \|^2} \)
		\State Set \( \alpha_t = \beta_t \)  
		\While{$f(x_t + \alpha_t d_t) > f(x_t) + \sigma \alpha_t ( g_t^T d_t - \alpha_t \mu L_t \| d_t \|^2 )$}
		\State Update \( \alpha_t \gets p \alpha_t \)
		\EndWhile
		\State \Return Step size \( \alpha_t \)
	\end{algorithmic}
\end{algorithm}
\begin{remark}
	In the modified Armijo line search, the parameter $L_t$ must be estimated at each iteration. For $t > 1$, we update $L_t$ as $L_t = \frac{s_{t-1}^T \bar{y}_{t-1}}{\|s_{t-1}\|^2}$.
	Comparing the classical Armijo line search with the modified Armijo line search (MALS), we observe that the step size $\alpha_t$ obtained from MALS \eqref{7777} ensures a greater descent magnitude in the objective function than that obtained from \eqref{7890}. The Modified Armijo line search (MALS)  balances step size efficiency and computational cost by avoiding excessively small step sizes for faster convergence, like Armijo, while eliminating the need for additional gradient evaluations, like Wolfe \cite{Xiao2013GlobalCO}, making it more computationally efficient.
\end{remark}
\section{M-LBFGS with displacement aggregation (AggMBFGS) and its convergence} \label{789}
AggMBFGS comprises two phases per iteration: (i) decision variable update and (ii) curvature pair update.\\

\textbf{Decision variable update:}
At iteration \( t \), we have access to the initial Hessian approximation \( G_0 \), the curvature pairs \((s_t, \bar{y}_t)\) (with at most \(\tau\) pairs), the current iterate \( x_t \), and the gradient \( \nabla f(x_t) \). Similar to L-BFGS, the inverse Hessian approximation is not computed explicitly. Instead, the descent direction \( d_t \) is determined at a computational cost of at most \( \mathcal{O}(\tau d) \) \cite{nocedal1999numerical}. The new iterate is then obtained using a step size \( \alpha_t \) as follows: 
\begin{equation}
	x_{t+1} = x_t - \alpha_t W_t \nabla f_t = x_t + \alpha_t d_t,
\end{equation}
where \( d_t = -W_t \nabla f_t \).   
Similar to L-BFGS, the memory requirement for M-LBFGS is \( \mathcal{O}(\tau d) \).  

\textbf{Curvature Pair Update:}  
The curvature pair update consists of three steps, which are described below.
\begin{enumerate}
	\item \textbf{Variable and gradient variation steps}:  we compute $s_t= x_{t+1}-x_{t}$, $ {y_{t}}=\nabla f(x_{t+1}) -\nabla f(x_{t})$, $\bar y_{t}=y_{t}+r_{t}\|g_{t}\|s_{t}$, $ r_{t}= 1+\max[0,-\frac {y_{t}^Ts_t}{s_t^T s_t}]$, and the descent direction $d_t$ using a two-loop recursion scheme described in Algorithm \ref{alg3}. The total computational complexity of Algorithm \ref{alg3} is $\mathcal{O}(4\tau d)$.
	
	\begin{algorithm}[ht]
		\caption{Computation of $d_t=-W_{t} \nabla f_{t}$ without explicitly computing $W_{t}$ \cite{nocedal1999numerical}}\label{alg3}
		\begin{algorithmic}
			\State{Initialize $g \leftarrow \nabla f_{t}$, $\rho_{t} \leftarrow \frac{1}{y_t^{T}s_t}, W_{t}^{0} \leftarrow$ initial inverse Hessian approximation}
			\For{$i=t-1, t-2, \ldots, t-\tau $}
			\State {$\alpha_{i} \leftarrow \rho_{i} s_{i}^{T} g$}
			\State    $g \leftarrow g-\alpha_{i} y_{i}$
			\EndFor
			\State $m \leftarrow W_{t}^{0} g $
			\For{$i=t-\tau, t-\tau+1, \ldots, t-1$}
			\State  $\beta \leftarrow \rho_{i} y_{i}^{T} m$
			\State   $m \leftarrow m+s_{i}\left(\alpha_{i}-\beta\right)$
			\EndFor
			\State stop with result $d_t=-m =- W_{t} \nabla f_{t}$.
		\end{algorithmic}
	\end{algorithm}
	\item \textbf{Displacement step}:
	Consider the set of curvature pairs \( P_{t-1} = \{(s_t, \bar{y}_t)\}_{t=0}^{\hat{\tau}-1} \), where $\hat{\tau}$ denotes the number of stored pairs, constrained by the memory size $\tau$. When a new curvature pair \((s_t, \bar{y}_t)\) is computed, the key question arises: How should \( P_t \) be updated to retain the most relevant information?  
	
	The displacement aggregation strategy integrates the new pair \((s_t, \bar{y}_t)\) into \( P_{t-1} \) to form \( P_t \) while removing dependent vectors. Rather than simply discarding the oldest pair \((s_0, \bar{y}_0)\), as in M-LBFGS \cite{Xiao2013GlobalCO} , this approach selectively incorporates \((s_t, \bar{y}_t)\). It has been shown that, under this strategy, the Hessian approximation in a limited-memory setting can be equivalent to a full-memory Hessian approximation~\cite{Berahas2019LimitedmemoryBW}. As a result, the impact of memory reduction diminishes, and the convergence rate improves. The implementation of displacement aggregation encompasses three cases:
	
	\noindent \textbf{Case 1} If the new variable variation $s_t$ is linearly independent of $\{s_i\}_{i=0}^{i=\hat \tau-1}$ in $P_{t-1}$, then a new curvature pair  $(s_t,\bar{y}_t)$ is added to $P_{t-1}$ to form $P_t$, i.e., 
	\begin{equation}\label{231}
		P_t=\{ P_{t-1}, (s_t,\bar{y}_t)\}.
	\end{equation}
	\textbf{Case 2} If the new variable variation  $s_t$ is linearly dependent on the previous variable variation $s_{\hat \tau-1}$ in $P_{t-1}$, then we update $P_t$ by removing the previous stored curvature pairs $(s_{\hat \tau-1}, \bar{y}_{\hat \tau-1})$ and replacing it with $(s_t,\bar{y}_t)$, i.e.,
	\begin{equation}
		P_t=\{(s_0, \bar{y}_0), \dots, (s_{\hat \tau-2}, \bar{y}_{\hat \tau-2}), (s_t, \bar{y}_t) \}.
	\end{equation}
	\textbf{Case 3}  If the new variable variation $s_t$ is linearly dependent on any previously stored variable variation in $P_{t-1}$, then we update $P_t$ by projecting the new pair $(s_t, \bar{y}_t)$ onto the subspace spanned by the existing pairs, thereby modifying $P_{t-1}$. Let us assume that $s_t=s_j$ where $0 \leq j \leq \hat \tau-1$.  Define
	\begin{equation*}
		S_{0:{\hat \tau}}=\{s_0, \dots, s_j, \dots, s_{\hat\tau}\}, \bar{Y}_{0:{\hat \tau}}=\{\bar{y}_0, \dots, \bar{y}_j, \dots, \bar{y}_{\hat \tau}\}  
	\end{equation*}
	and let $(s_{\hat \tau}, \bar{y}_{\hat \tau})=(s_t, \bar{y}_t)$.
	Assume that \((s_{\hat{\tau}}, \bar{y}_{\hat{\tau}})\) represents the new curvature pairs. Then, to update the set of curvature pairs, we remove the pair \((s_j, \bar{y}_j)\) and replace the subsequence \(\bar{Y}_{j+1:\hat{\tau}}\) with its modified version \(\hat{Y}_{j+1:\hat{\tau}}\). Consequently, the updated set \(P_t\) is expressed as  
	\begin{equation}\label{666}
		P_t = \{(s_0, \bar{y}_0), \dots, (s_{j-1}, \bar{y}_{j-1}), (s_{j+1}, \hat{y}_{j+1}), \dots, (s_{\hat{\tau}}, \hat{y}_{\hat{\tau}}) \}.
	\end{equation}
	The computation of \(\hat{Y}_{j+1:\hat{\tau}}\) follows the method presented in \cite{Berahas2019LimitedmemoryBW}, and  is given by  
	\begin{equation}
		\hat{Y}_{j+1:\hat{\tau}} = (W_{0:j-1})^{-1} S_{j+1:\hat{\tau}}
		\begin{bmatrix}
			A & 0
		\end{bmatrix}
		+ \bar{y}_j 
		\begin{bmatrix}
			b \\ 
			0
		\end{bmatrix}^{T}
		+ \bar{Y}_{j+1:\hat{\tau}},
	\end{equation}
	where \(W_{0:j-1}\) represents the inverse Hessian approximation computed using the limited-memory approach, initialized with \(W_0\) and updated with \(j\) curvature pairs \((s_i, \bar{y}_i)_{i=0}^{i=j-1}\). The matrix \(A \in \mathbb{R}^{(\hat{\tau}-j) \times (\hat{\tau}-j-1)}\) and the vector \(b \in \mathbb{R}^{\hat{\tau}-j-1}\) can be determined using Algorithm 4 and Algorithm 3 from \cite{Berahas2019LimitedmemoryBW}, respectively. Additionally, the term \((W_{0:j-1})^{-1} S_{j+1:\hat{\tau}}\) is also derived as per \cite{Berahas2019LimitedmemoryBW}. To identify the appropriate case during the update, the Cholesky factorization of the inner product matrix corresponding to the stored displacement vectors \(\{(s_0, \bar{y}_0), \dots, (s_{\hat{\tau}-2}, \bar{y}_{\hat{\tau}-2}), (s_{\hat{\tau}-1}, \bar{y}_{\hat{\tau}-1})\}\) can be computed. Implementation details are thoroughly discussed in \cite{Berahas2019LimitedmemoryBW}. Case-1 occurs only if \(\hat{\tau} \leq \tau\). If $\hat{\tau} \geq \tau$, then $s_t$ cannot be linearly independent of $\{s_i\}_{i=0}^{i=\hat{\tau}-1}$ in $P_{t-1}$ because all the vectors are selected from the subset \(\{e_1, \dots, e_{\tau}\}\), which has size $\tau$ . Hence, the number of curvature pairs in $P_t$ remains bounded by the memory size \(\hat{\tau} +1 \leq \tau\). In both Case-2 and Case-3, the number of curvature pairs remains $\hat{\tau}$ and does not increase to $\hat{\tau}+1$. Hence, the number of curvature pairs in $P_t$ is $\hat{\tau}$, which is bounded by the memory size $\tau$. The computational complexities for these cases are summarized in Table \ref{tab4}.
	\begin{table}[htbp]
		\caption{Computational complexity to compute A, b, $(W_{0:j-1})^{-1} S_{j+1:{\hat \tau}}$, Cholesky factorization of an inner product matrix for checking Case-1, Case-2, and Case-3,  $\{G_t e_i\}_{i=1}^{i=\tau}$ , $d_t$}\label{tab4}
		\centering
		\begin{tabular}{|c|c|}
			\hline
			Parameters & Computational complexity  \\
			\hline
			A & $\mathcal{O}(\tau^2d+\tau^4)$ \cite{Berahas2019LimitedmemoryBW}\\ 
			\hline
			b & $\mathcal{O}(\tau^2 d)$ \cite{Berahas2019LimitedmemoryBW} \\
			\hline
			$(W_{0:j-1})^{-1} S_{j+1:{\hat \tau}}$&  $\mathcal{O}(\tau^2 d)$ \cite{Berahas2019LimitedmemoryBW} \\
			\hline
			Cholesky factorization of an inner
			product matrix 
			for checking 3 cases &  $\mathcal{O}(\tau d)$ \cite{Berahas2019LimitedmemoryBW}\\
			\hline
			$\{G_t e_i\}_{i=1}^{i=\tau}$ & $\mathcal{O}(\tau^2 d)$ (Check section 7.2 in \cite{nocedal1999numerical})\\
			\hline
			$d_t$ & $\mathcal{O}(\tau d)$ \cite{nocedal1999numerical}\\ 
			\hline
		\end{tabular}
	\end{table}
\end{enumerate}
\begin{algorithm} 
	\caption{AggMBFGS with Backtracking Modified Armijo Line Search} \label{1515}
	\begin{algorithmic}[1]
		\Require $x_0$, initial inverse Hessian approximation $W_0$, line search parameters $p \in (0,1)$, $\mu \in [0,\infty)$, $\sigma \in (0,1)$
		\For{$t=0,1,\dots,T$}
		\State Compute the gradient $\nabla f(x_t)$
		\State Update $x_{t+1} = x_t +\alpha_{t} d_t$, where $d_t$ is computed by Algorithm~\ref{alg3}
		\State Compute step length $\alpha_t>0$ using Algorithm~\ref{alg:2}
		\State Compute curvature pairs $(s_t,\bar{y}_t)_{t=0}^{\hat{\tau}-1}$: $s_t = x_{t+1}-x_{t}$  and  $\bar{y}_t = y_t + r_t \|g_t\|s_t$ where $y_{t}= \nabla f(x_{t+1}) - \nabla f(x_{t})$ and $r_{t}=1+ \max\left(0,-\frac {y_{t}^Ts_t}{s_t^T s_t}\right)$
		\State Update historical curvature pair $P_t$ using equations~(\ref{231})--(\ref{666})
		\EndFor
	\end{algorithmic}
\end{algorithm}
\subsection{Convergence analysis of AggMBFGS } \label{1234}
\begin{assumption}
	\label{39999}
	Let us assume that level set is $ \Omega =\{x \mid f(x)\leq f(x_0) \}$ .
	\begin{enumerate}
		\item  $f(x)$  is  bounded below and the level set $\Omega $ is bounded.
		
		\item  $g(x)=\nabla f(x)$ is Lipschitz continuous in an open neighborhood $N$ of the level set $\Omega$ i.e.,
		\begin{equation*}
			\|g(x)-g(y)\|< L \|x-y\| \quad \forall x,y\in N
		\end{equation*}
		where L is Lipschitz constant.
	\end{enumerate}
\end{assumption}
\subsection*{Local superlinear convergence of MBFGS}
\begin{proposition}\label{10101}
	Suppose that the level set $\Omega = \{ x : f(x)\leq f(x_0) \} $ is bounded and  $ f(x) $ is twice continuously differentiable near $x^*$, which is contained in $\Omega$. Let $x_{t} \rightarrow  x^* $ where $g(x^*)=0$, Hessian $G(x)$ is positive definite and Lipschitz continuous at $x^*$ and $\alpha_{t} \in (0,1/2) $ is satisfied by backtracking modified Armijo line search \ref{alg:2}, then the sequence $\{x_{t}\}$ generated by MBFGS converges to $x^*$ superlinearly.
\end{proposition}
\begin{proof}    
	The proof follows the same structure as the proof of Theorem 3.8 in \cite{Li2001AMB}.
\end{proof}	
\begin{remark}
	Such a result cannot be proved for the M-LBFGS. The modified limited memory BFGS method achieves a local linear rate, which is not better than the gradient descent method. One can refer to \cite{Nocedal1980UpdatingQM} for the proof.
\end{remark}
\subsection*{Global convergence of MBFGS}
\begin{proposition}\label{8978}
	Let $\{x_t\}$ be generated by MBFGS with backtracking modified Armijo line search  \ref{alg:2} and satisfy Assumption \ref{39999}. Then we have 
	\begin{equation}\label{400}
		\liminf\limits_{k\rightarrow \infty} \|g_t\|=0,
	\end{equation} 
	i.e., there exists a point $x^*$ such that $x_t \rightarrow x^*$.
\end{proposition}
\begin{proof}
	One can refer to the proof of Theorem 5.1 in \cite{Li2001AMB}.
\end{proof}
\begin{theorem}\label{3555}
	If AggMBFGS and MBFGS are performed with a memory size $\tau$, and have the same initial settings, then the iterates generated by AggMBFGS are equal to those generated by MBFGS.
\end{theorem}
\begin{proof}
 This result follows from the proof of Lemma 1 in \cite{Gao2023LimitedMemoryGQ}. Specifically, replacing the curvature pair $(s_t,r_t)$ with $(s_t,\bar{y}_t)$ in the proof of Lemma 1 in \cite{Gao2023LimitedMemoryGQ} establishes the desired conclusion.
\end{proof}
Using Propositions~\ref{10101} and \ref{8978}, along with Theorem~\ref{3555}, we establish the local superlinear rate of convergence and global convergence of AggMBFGS with a backtracking MALS for a sufficiently large memory size \( \tau \), assuming the same initial settings as MBFGS.
\subsection*{Global Convergence of M-LBFGS}
The global convergence of M-LBFGS with the Wolfe line search was studied in \cite{Xiao2013GlobalCO}, while its convergence with the Armijo line search was analyzed in \cite{Shi2016ALM}. Since the Modified Armijo line search accelerates convergence by preventing excessively small step sizes, like the Armijo line search, while preserving computational efficiency by avoiding extra gradient evaluations required in the Wolfe line search, we aim to establish the global convergence of M-LBFGS with the Modified Armijo line search.
Let $\theta_{t}$ be the angle between $s_t$ and $G_t s_t$, then 
$cos(\theta_{t})=\frac{s_t^T G_t s_t}{\|s_t\| \|G_t s_t \|}= - \frac{g_t^T s_t}{\|s_t\| \|g_t\|}$ \cite{Shi2016ALM}.
\begin{lemma}
	\label{9999}
	Let us choose $G_0$ in such a way that $\|G_0\|$ and $\|G_0^{-1}\|$ are bounded and $\{x_t\}$ be the sequence of iterate  generated by Algorithm \ref{506}. if $\|g_t\| \geq \epsilon$ holds $\forall t$ with some non-negative constant $\epsilon \geq 0$, then their exist positive constant $q>0$ such that the inequality 
	\begin{equation}
		cos(\theta_{t}) \geq q
	\end{equation}
	hold $\forall t$.
\end{lemma}
\begin{proof}
	One can refer to the proof of Lemma 3.3 in \cite{Xiao2013GlobalCO}.
\end{proof}
\begin{theorem}
	Let $\{x_t\}$ be generated by Algorithm-\ref{506} and satisfy  Assumption-\ref{39999}. Then we have 
	\begin{equation}\label{403}
		\liminf\limits_{t\rightarrow \infty} \|g_t\|=0,
	\end{equation} 
	i.e., there exists a point $x^*$ such that $x_t \rightarrow x^*$.
\end{theorem}
\begin{proof}
	On the contrary, assume that \eqref{403}  doesn't hold, i.e., there is a constant $\epsilon >0$ such that $\|g_t\|\geq \epsilon~\forall t $.  If $\alpha_{t} \neq \beta_t$, it follows from step 4 of Algorithm \ref{1515} that $p^{-1}\alpha_{t}$ doesn't satisfy  modified Armijo, i.e.,
	\begin{equation}\label{1556}
		f(x_t + p^{-1}\alpha_{t}d_t) - f(x_t) > p^{-1} \sigma \alpha_{t}(g_t^{T}d_t -p^{-1} \alpha_{t} \mu L_t \| d_t\|^2 ).
	\end{equation} 
	By mean value theorem, there exist $\theta_{t} \in (0,1)$ such that 
	\begin{eqnarray*}
		f(x_t + p^{-1}\alpha_{t}d_t) - f(x_t)&=& p^{-1} \alpha_{t} g(x_t+\theta_{t}p^{-1}\alpha_{t}d_t)^T d_t\\
		&=& p^{-1} \alpha_{t} g_t^{T}d_{t}+p^{-1}\alpha_{t}( g(x_t+\theta_{t}p^{-1}\alpha_{t}d_t) -g(x_t))^T d_t\\
		&\leq& p^{-1} \alpha_{t} g_t^{T}d_{t} + Lp^{-2} \alpha_{t}^{2}\|d_t \|^2.
	\end{eqnarray*}
	Using Equation (\ref{1556}), we have 
	\begin{eqnarray*}
		\sigma (g_t^{T}d_t -p^{-1} \alpha_{t} \mu L_t \| d_t \|^2 ) 
		&<&  g_t^{T}d_{t} + Lp^{-1} \alpha_{t} \|d_t\|^2  \\ 
		\Rightarrow \alpha_{t} & >& \frac {p(\sigma -1)g_t^{T}d_t}{(L+\sigma\mu L_t) \|d_t\|^2}.
	\end{eqnarray*}
	So, we know that from (\ref{7777})  
	\begin{eqnarray*}
		f(x_t + \alpha_{t}d_t) - f(x_t) &<& \sigma \frac {p(\sigma -1)(g_t^{T}d_t)^2}{(L+\sigma\mu L_t) \|d_t\|^2} -\frac{\sigma \mu L_t (1-\sigma)^2 p^2 (g_t^{T}d_t)^2}{(L+\sigma\mu L_t)^2 \|d_t\|^2}\\
		\Rightarrow f(x_t)-f(x_t + \alpha_{t}d_t) &>& \frac {p\sigma(1-\sigma)(g_t^{T}d_t)^2}{(L+\sigma\mu L_t) \|d_t\|^2}+\frac{\sigma \mu L_t (1-\sigma)^2 p^2 (g_t^{T}d_t)^2}{(L+\sigma\mu L_t)^2 \|d_t\|^2}\\
		&>&M \|g_t\|^2 cos^{2}\theta_{t}, 
	\end{eqnarray*}
	where $M= \frac{p\sigma(1-\sigma)}{(L+\sigma\mu L_t)}+ \frac{\sigma \mu L_t (1-\sigma)^2 p^2}{(L+\sigma\mu L_t)^2}$.
	Then $\forall t>0$. We have
	\begin{equation}
		\label{500}
		\sum\limits_{i=0}^{t-1}  f(x_i)-f(x_i + \alpha_{i}d_i)=f(x_0)-f(x_t) > M \sum\limits_{i=0}^{t-1} \|g_t\|^2 cos^{2}\theta_{t}. 
	\end{equation}
	Taking limit on both sides as $t \rightarrow \infty$ and $f(x_t)\to f(x^*)$, we get $\liminf\limits_{t\rightarrow \infty} \|g_t\| cos\theta_{t}=0$. 
	From Lemma \ref{9999}, we deduce  $\liminf\limits_{t\rightarrow \infty} \|g_t\|=0$. This contradicts our assumption. Similarly, if $\alpha_{t} =\beta_t$ where $\beta_t= -g_{t}^T d_t/(L_{t}\|d_t \|^2)$, then it follows from \eqref{7777} that
	\begin{eqnarray*}
		f(x_t + \alpha_{t}d_t)&\leq& f(x_t)+\sigma \beta_t (g_t^{T}d_t -\beta_{t} \mu L_t \| d_t \|^2 ) \\
		&=&f(x_t)+\sigma (\mu+1) \beta_t g_t^{T}d_t \\
		& =&  f(x_t)- \frac{\sigma (\mu+1)} { L_t \| d_t\|^2} \|g_t\|^2 \|d_t\|^2 cos^2{\theta_{t}}.
	\end{eqnarray*}
	Then $ f(x_t) -f(x_t + \alpha_{t}d_t) \geq \frac{\sigma (\mu+1)} { \|L_t\|} \|g_t\|^2 cos^2{\theta_{t}} = N \|g_t\|^2 cos^2{\theta_{t}}$,
	where $N=\frac{\sigma (\mu+1)} { \|L_t\|}$.
	Similarly, $\forall t>0$
	\begin{equation}
		\label{5555}
		\sum\limits_{i=0}^{t-1}  f(x_i)-f(x_i + \alpha_{i}d_i)=f(x_0)-f(x_t) > N \sum\limits_{i=0}^{t-1} \|g_t\|^2 cos^{2}\theta_{t}. 
	\end{equation}
	Taking limit on both sides as $t \rightarrow \infty$ and $f(x_t)\to f(x^*)$, we get $\liminf\limits_{t\rightarrow \infty} \|g_t\| cos\theta_{t}=0$. 
	From Lemma \ref{9999}, we conclude  $\liminf\limits_{t\rightarrow \infty} \|g_t\|=0$. This contradicts our assumption. Hence, considering both cases, we get our desired results.
	
\end{proof}

\subsection{\textbf{Computational complexity of AggMBFGS}}  
The total computational cost of the M-LBFGS method is \( \mathcal{O}(5 \tau d) \) per iteration \cite{Liu1989OnTL}, where \( d \) is the number of variables in the optimization algorithm, and \( \tau \) is the user-defined memory allocation parameter, typically chosen in the range \( 3 \leq \tau \leq 10 \). AggMBFGS incorporates displacement steps, which are straightforward to implement. The total computational cost of the modified \( \hat{Y}_{j+1:{\hat \tau}} \) is $
\mathcal{O}(\tau^{2} d + \tau^{4})$  
as analyzed in \cite{Berahas2019LimitedmemoryBW}. 

\begin{table}[htbp]
	\caption{Comparison between AggMBFGS and other quasi-newton methods}
	\centering
	\begin{tabular}{|c|c|c|c|}
		\hline
		Algorithm & Types of function & Memory & Complexity\\
		\hline
		AggMBFGS(this work) & $C^2$+ non-convex& $\mathcal{O}(\tau d)$ & $\mathcal{O}(\tau^{2} d+ \tau^{4} )$ \\ 
		\hline
		L-BFGS with DA \cite{Berahas2019LimitedmemoryBW}& $C^2$ +uniformly convex function& $\mathcal{O}(\tau d)$ &  $\mathcal{O}(\tau^{2} d+ \tau^{4} )$ \\ 
		\hline
		L-BFGS \cite{Liu1989OnTL} & $C^2$+uniformly convex function  & $\mathcal{O}(\tau d)$  & $\mathcal{O}(\tau d)$ \\
		\hline
		Modified BFGS \cite{Li2001AMB} & $C^2$+ non-convex  & $\mathcal{O}( d^2 )$  &  $\mathcal{O}( d^2 )$  \\ 
		\hline
		Classical BFGS \cite{nocedal1999numerical} & $C^2$+uniformly convex function & $\mathcal{O}(d^2)$ &  $\mathcal{O}( d^2 )$   \\ 
		\hline
		
	\end{tabular}
\end{table}

\subsection{Comparison between AggMBFGS and M-LBFGS} \label{777}  

In this section, we examine the effectiveness of Algorithm \ref{506}  and Algorithm \ref{1515}. A collection of 47 nonlinear unconstrained problems is used in our experiment. We conducted numerical experiments using the CUTEst environment \cite{gould2015cutest}, with all test problems obtained from CUTEst.
\newpage
 Throughout the paper, we use the following computational setup:  \\
\textbf{System}: Intel Core i5-10210U, 2.11 GHz, 8GB RAM, Ubuntu Linux. \\ 
 \textbf{Interface}: MATLAB 2020b.\\
  \textbf{Initialization}  $
	x_0 = \frac{\text{randn}(d,1)}{\lVert \text{randn}(d,1) \rVert}
	$, $G_0 = \mathbf{I}_{d \times d}$.\\
\textbf{Parameters}: CPU time limit = 3600s, $\mu = 1$, $\sigma = 0.2$, $p = 0.3$, $L_0 = 1$.  \\
 \textbf{Termination criteria}: $\|g_t\|_\infty \leq 10^{-6} \max(1,\|g_0\|_\infty)$ or iteration count $\geq 10^5$.  
 \par We have taken memory size $\tau=5$. We tested problems with dimensions ranging from 2 to 132,200. In the majority of cases, \( d \gg \tau \), which means that the computational cost of our aggregation scheme is negligible compared to the computational cost of calculating search directions. The comparison between AggMBFGS and M-LBFGS is presented in Table \ref{909086} and Table \ref{7873}. We observe from Table \ref{909086} and Table \ref{7873} that Algorithm \ref{1515} outperforms Algorithm \ref{506} in the number of iterations and function evaluations.
\begin{table}[htbp]
	\centering
	\footnotesize
	\caption{Name of the problems, dimension of the problem, (number of iteration, function evaluation, and aggregation used in AggMBFGS when applied to solve the problem from CUTEst set with $n \in [2,123200]$),  (number of iteration, function evaluation in M-LBFGS when applied to solve the problem from CUTEst set with $n \in [2,123200]$)}.
	\label{909086}
	\begin{tabular}{l*{5}{l}{l}}\toprule
		\multirow{2}{*}{Name}  &  \multirow{2}{*}{Dim} & \multicolumn{3}{l}{ AggMBFGS[5]  } & \multicolumn{2}{l}{ M-LBFGS[5] }\\
		\cmidrule(lr){3-5} \cmidrule(lr){6-7}
		& (n) & Iters & func & Agg & Iters & Func\\	
		\midrule
		ARGLINA &200 &2 &3 &0 &2 &3 \\ \\
		ARGLINB &200 &3 &49&0 &3&49\\ \\
		ARGLINC &200 &3 &49&0 &3&49\\ \\
		ARGTRIGLS&200& 669& 10362 & 0& 669& 10362\\ \\
		ARWHEAD &5000&49&212&43&59&235\\ \\
		BA-L1LS&57&352&8284&0&352&8284\\ \\
		BA-L1SPLS&57&50&1549&0&50&1549\\ \\		BDQRTIC&5000&62&718&0&62&718\\
		\\		BOX &10000&207&1399&106&280&1930\\ \\
		BOXPOWER&20000&10&71&1&10&71\\ \\
		BROWNAL&200&17&114&1&17&115\\ \\
		BROYDN3DLS&5000&89&429&0&89&429 \\ \\
		BROYDN7D&5000&7830&49931&0&7830&49931\\ \\
		BROYDNBDLS&5000&104&756&0&104&756\\ \\
		BRYBND&5000&104&756&0&104&756\\ \\
		CHAINWOO&4000&353&3630&0&353&3630\\ \\
		CHNROSNB&50&335&3065&0&335&3065\\ \\
		CHNRSNBM&50&344&3200&0&344&3200\\ \\
		CURLY10&10000&2183&28121&0&2183&28182\\ \\
		CURLY20&10000&417&6167&0&417&6167\\ \\
		CURLY30&10000&326&5383&0&326&5383\\ \\
		DIXMAANA&3000&15&40&9&19&52\\ \\
		DIXMAANB&3000&32&63&0&32&63\\ \\

		\bottomrule
	\end{tabular}
\end{table}
\begin{table}[htbp]
	\centering
	\footnotesize
	\caption{Name of the problems, dimension of the problem, (number of iteration, function evaluation, and aggregation used in AggMBFGS when applied to solve the problem from CUTEst set with $n \in [2,123200]$),  (number of iteration, function evaluation in M-LBFGS when applied to solve the problem from CUTEst set with $n \in [2,123200]$.)}
	\label{7873}
	\begin{tabular}{l*{5}{l}{l}}\toprule
		\multirow{2}{*}{Name}  &  \multirow{2}{*}{Dim} & \multicolumn{3}{l}{ AggMBFGS[5]  } & \multicolumn{2}{l}{ M-LBFGS[5] }\\
		\cmidrule(lr){3-5} \cmidrule(lr){6-7}
		& (n) & Iters & func & Agg & Iters & Func\\	
		\midrule
		DIXMAANC&3000&25&57&0&25&57\\ \\
		DIXMAAND&3000&19&50&0&19&50\\ \\
		DIXMAANE&3000&268&498&0&268&498\\ \\
		DIXMAANF&3000&203&250&0&203&250\\ \\
		DIXMAANG&3000&132&224&0&132&224\\ \\
		DIXMAANH&3000&288&1066&0&288&1066\\ \\
		DIXMAANI&3000&201&230&0&201&230\\\\
		DIXMAANJ&3000&128&253&0&128&253\\ \\ 
		DIXMAANK&3000&132&301&0&132&301\\ \\
		DMN15333LS&99&12&253&1&13&253\\ \\
		DMN37142LS&66&8531&80501&0&8531&80501\\ \\
		ERRINROS&50&60&604&0&60&604\\ \\
		ERRINRSM&50&74&721&0&74&721\\ \\
		FLETBV3M&5000&10&86&3&12&141\\ \\
		HILBERTA&2&13&48&10&15&50\\ \\
		NONCVXU2&5000&34&138&6&36&142\\ \\
		NONDQUAR&5000&33&231&0&33&231\\ \\
		POWELLSG&5000&28&269&38&39&155\\ \\
		POWER&10000&55&1841&0&55&1841\\ \\
		SPARSQVR&10000&78&1018&0&78&1018\\ \\
		TESTQUAD&5000&4866&95084&0&4866&95084 \\ \\
		TQUARTIC&5000&236&2939&633&342&1434\\ \\
		YATP2LS&123200&264&3062&0&264&3062\\ \\
		YATP1LS&123200&46&400&25&56&294\\
		\bottomrule
	\end{tabular}
\end{table}

\subsection{Application of AggMBFGS in finding the largest Eigenvalue problem}\label{478}


Here, we focus on applying the AggMBFGS algorithm with a modified Armijo line search to minimize  
\[
f(x) = \frac{\|x \|^4}{4}- \frac { x^{T}Ax}{2}
\]  
for computing the largest eigenvalue of real positive definite matrices. The gradient and Hessian of \( f(x) \) are given by:  

\begin{equation*}
	\nabla f(x) = \|x\|^2 x - Ax,
\end{equation*}
and  
\begin{equation*}
	\nabla^2 f(x) = \|x\|^2 \mathbb{I}_d + 2 x x^T - A,
\end{equation*}  
where \( \mathbb{I}_d \) is the identity matrix of order \( d \).  

\par It is straightforward to verify that \( f(x) \) is bounded below, the level set  
\[
\Omega = \{x \mid f(x) \leq f(x_0) \}
\]  
is bounded, and \( \nabla f(x) \) is Lipschitz continuous in an open neighborhood of \( \Omega \) \cite{Shi2016ALM}. Hence, \( f(x) \) satisfies Assumption~\ref{39999}. Therefore, we can apply Algorithm \ref{1515} to minimize \( f(x) \).  

\par Moreover, from \cite{Mongeau2004ComputingEO}, Theorem~2.3, if \( x^* \) is the global minimum and the corresponding critical value satisfies  
\[
f(x^*) = -\frac{\lambda^2}{4},
\]  
then \( \lambda \) is the largest eigenvalue of \( A \).
\newpage
\section{Numerical experiment}\label{3456}

The main objective of this experiment is to demonstrate the effectiveness of AggMBFGS in computing the largest eigenvalue of real symmetric positive definite matrices. This study compares AggMBFGS with the \texttt{eigs} MATLAB command. We select 37 large-scale real symmetric positive definite sparse matrices, ranging in size from 4,098 to 54,929, from \url{http://www.cise.ufl.edu/research/sparse/matrices} \cite{davis2011university} and evaluate them using Algorithm \ref{1515}. We estimate the largest eigenvalue of positive definite matrices using the \texttt{eigs} MATLAB command in advance to compute the relative error, defined as:  
\begin{equation*}
	\text{Relative error} = \frac{\lVert{\texttt{eigs}}(\lambda) - \lambda \rVert}{\lambda},
\end{equation*}  
where \( \lambda \) is the largest eigenvalue computed by AggMBFGS, and \( {\texttt{eigs}}(\lambda) \) is the largest eigenvalue obtained using the \texttt{eigs} MATLAB command. 

\par In Table \ref{tab:1} and Table \ref{tab:2}, we compare the largest eigenvalues computed by Algorithm \ref{1515} and M-LBFGS with classical Armijo line search \cite{Shi2016ALM} with those obtained using the \texttt{eigs} MATLAB command. The results demonstrate that Algorithm \ref{1515} outperforms M-LBFGS \cite{Shi2016ALM} regarding accuracy with respect to the \texttt{eigs} command. Therefore, we recommend applying Algorithm \ref{1515} for computing the largest eigenvalue of large real symmetric positive definite matrices. Algorithm \ref{1515} is more reliable due to its lower storage requirements and faster convergence rate, making it well-suited for handling large real symmetric positive definite matrices.

\begin{remark}
In \cite{Shi2016ALM}, $\tau=3$ outperforms other values of $\tau$ in the implementation of M-LBFGS. Therefore, we set $\tau=3$ while implementing AggMBFGS with a backtracking modified Armijo line search to evaluate the largest eigenvalue of real positive definite matrices.

\end{remark}


\begin{table}[htbp]
	\centering
	\footnotesize
	\caption{Name of the matrix, Order of the matrix, Number of iteration, Time taken in second, The largest eigenvalue using AggMBFGS with modified Armijo line search, Largest eigenvalue using the \texttt{eigs} Matlab command, \bf{Relative error} \normalfont {of the largest eigenvalue obtained between \texttt{eigs} and }\bf{AggMBFGS[3] using modified Armijo line search}, \bf{Relative error} \normalfont {of the largest eigenvalue obtained between \texttt{eigs} and} \bf{ M-LBFGS [3] using Armijo line search}.}
	\label{tab:1}
	\begin{tabular}{l*{5}{l}{l}}\toprule
		\multirow{2}{*}{Name}  &  \multirow{2}{*}{Order} & \multicolumn{2}{l}{ AggMBFGS[3]  } & \multicolumn{3}{l}{ \texttt{eigs}(MATLAB) }{ { Rel. error with respect to \texttt{eigs} }}\\
		\cmidrule(lr){3-4} \cmidrule(lr){5-5} \cmidrule(lr){6-7}
		& (d) & Iter/Time &  largest eigenvalue &  largest eigenvalue & AggMBFGS[3] &  M-LBFGS[3] \\	
		\midrule
		c-30 &5321 &5/0.05 &4.68199102e+6 & 4.68199099e+6&6.4075e-9& 1.3353e-7\\ \\
	c-33 &6317 &6/0.09 &2.04540851e+5 &2.04540830e+5&1.0267e-7& 5.8454e-7 \\ \\
	c-36&7479& 9/0.35 & 8.74111892e+3 & 8.74111803e+3 & 1.0182e-7& 1.7819e-7\\ \\
	c-50 &22401&8/0.26&2.65859662e+5&2.65859669e+5&2.6330e-8& 7.5366e-7\\ \\
	bloweybl & 30003 &8/0.19 &1.00007499e+2 & 1.00007501e+2& 1.9999e-8& 5.3367e-7 \\ \\
	net150 & 43520 & 37/10.24 &1.45592934e+2 &1.45593064e+2& 8.9290e-7&9.6467e-7  \\  \\
	mark3jac100sc & 45769&12/0.65 &1.04857598e+6 &1.04857600e+6& 1.9073e-8& 4.0455e-9 \\  \\
	rajat27&20640& 19/0.35 &7.69115075e+5&7.69114987e+5&1.1442e-7& 7.3177e-7\\ \\
	rajat01&6833& 11/1.05& 4.21268417e+1&4.21268444e+1& 6.4092e-8& 8.5859e-7\\ \\
	sts4098&4098& 15/0.16 &3.07102523e+8&3.07102503e+8&  6.5125e-8& 6.3120e-7 \\ \\
	bcsstk28& 4410 & 9/0.22&7.69621465e+8&7.69621402e+8& 8.1858e-8&5.5987e-7  \\ \\
	mhd4800b& 4800 & 7/0.14 &2.19626865e+0&2.19626863e+0&9.1064e-9& 8.2787e-7  \\ \\
	bcsstk16& 4884 & 27/1.04&4.94316578e+9&4.94316563e+9& 3.0345e-8& 8.8428e-7  \\ \\
	bloweybq & 10001& 8/0.10 &4.99974999e+3 &4.99975005e+3&1.2001e-8  &7.3897e-7   \\ \\
	ecl32&51993&29/4.45&9.61855290e+3&9.61854782e+3&5.2815e-7& 6.9659e-7\\ \\
	as-22july06 & 22963& 32/1.41 &7.1612987e+1&7.16130003e+1&1.8572e-8    &8.9177e-7  \\   \\ 
	ca-CondMat & 23133 & 29/1.81  &3.79541144e+1 &3.79541129e+1&3.9521e-8  &6.2849e-7  \\ \\
	mult-dcop-02 &25187 & 8/0.36 &1.25613773e+3 &1.25613771e+3& 1.5922e-8  & 8.5984e-8  \\     \\ 
		c-52&23948&6/0.18&1.94346250e+15&1.94346249e+15&5.1455e-9 &7.3524e-7\\ \\
	\end{tabular}
\end{table}
\begin{table}[htbp]
	\centering
	\footnotesize
	\caption{Name of the matrix, Order of the matrix, Number of iteration, Time taken in second, The largest eigenvalue using AggMBFGS with modified Armijo line search, Largest eigenvalue using the \texttt{eigs} Matlab command, \bf{Relative error} \normalfont {of the largest eigenvalue obtained between \texttt{eigs} and }\bf{AggMBFGS[3] using modified Armijo line search}, \bf{Relative error} \normalfont {of the largest eigenvalue obtained between \texttt{eigs} and} \bf{ M-LBFGS [3] using Armijo line search}.}
	\label{tab:2}
	\begin{tabular}{l*{5}{l}{l}}\toprule
		\multirow{2}{*}{Name}  &  \multirow{2}{*}{Order} & \multicolumn{2}{l}{ AggMBFGS[3]  } & \multicolumn{3}{l}{ \texttt{eigs}(MATLAB) }{ { Rel. error with respect to \texttt{eigs} }}\\
		\cmidrule(lr){3-4} \cmidrule(lr){5-5} \cmidrule(lr){6-7}
		& (d) & Iter/Time &  largest eigenvalue &  largest eigenvalue & AggMBFGS[3] &  M-LBFGS[3] \\	
		\midrule
			c-53&30235&12/0.90 &5.36188735e+03&5.36188723e+03&2.2380e-8& 5.3364e-7\\ \\
		c-54&31793&4/0.20&1.82553378e+08&1.82553380e+08&1.0965e-8&1.9374e-7 \\ \\
		c-56 &35910&6/0.31&1.20695037e+05&1.20695045e+05&6.6283e-8 &4.0152e-7 \\ \\
		c-57&37833&10/0.91&7.27009766e+04&7.27009867e+04&1.3893e-7& 2.6073e-7\\ \\
		c-58&37595&6/0.41&5.41547153e+04&5.41547139e+04&2.5852e-8&7.3685e-7\\ \\
		c-59&41282&20/1.59&8.38556378e+03&8.38556378e+03 &8.0012e-10 &4.5613e-7\\ \\
		c-65&48066&19/1.54&1.31413054e+05&1.31413052e+05&1.5219e-8& 2.1780e-8 \\ \\
		c-66&49989&14/1.23&1.70373624e+04&1.70373629e+04&2.9347e-8 & 5.2389e-7\\ \\
		c-64b&51035&7/0.94&2.00080883e+05&2.00080850e+05&2.9347e-8 & 3.4414e-7\\ \\
		bcsstk17 &10974 & 45/2.55 &1.29606157e+10 &1.29606158e+10&7.7157e-9 &8.3157e-7  \\ \\
		bcsstk25 & 15439 & 17/0.74 &1.06002059e+15&1.06002050e+15&8.4904e-8 &9.9397e-7 \\ \\
		olafu & 16146& 6/0.54 &9.47870244e+11 &9.47870339e+11&1.0022e-7  & 4.1015e-7  \\ \\
		gyro-k & 17361& 30/4.41 &3.65695228e+09&3.65695233e+09& 1.3673e-8& 7.5972e-7     \\   \\ 
		gyro & 17361&30/4.34 &3.65695228e+09 &3.65695233e+09&1.3673e-8& 7.5974e-7 \\ \\
		rajat26&51032&15/1.11 &8.26006295e+05&8.26006365e+05&8.8745e-8&2.9696e-7\\ \\
		rajat22&39899&21/1.14 &9.98382892e+05&9.98382851e+05&4.1066e-8& 9.3670e-7\\ \\
		rajat15&37261&16/1.12&3.26846630e+05&3.26846609e+05&6.4250e-8& 8.2989e-7\\ \\
		net100 & 29920 & 32/7.34 & 1.22667844e+02 &1.22667848e+02&3.2608e-8&9.1562e-7 
	\end{tabular}
\end{table}

\subsection{Error analysis}

From Table \ref{tab:1} and Table \ref{tab:2}, we observe that the largest eigenvalue computed using Algorithm \ref{1515} is close to that obtained with the \texttt{eigs} command in MATLAB. In \cite{Shi2016ALM}, Zhanwen Shi et al. demonstrated that the largest eigenvalue can be computed using M-LBFGS with the classical Armijo line search. Our proposed Algorithm \ref{1515} achieves a smaller relative error than M-LBFGS with the classical Armijo line search \cite{Shi2016ALM}, further validating its accuracy and effectiveness.

\section{Conclusion} \label{5876}

We have established the global convergence of both the AggMBFGS and M-LBFGS methods, as well as the local superlinear rate of convergence of AggMBFGS when used with the backtracking modified Armijo line search. Algorithm \ref{1515} (AggMBFGS with backtracking modified Armijo line search) demonstrates superior performance compared to M-LBFGS with the backtracking modified Armijo line search when applied to test problems from the CUTEst environment. Additionally, we have successfully employed Algorithm \ref{1515} to compute the largest eigenvalue of high-dimensional real positive definite matrices.  

\par The relative errors are computed with respect to the \texttt{eigs} MATLAB command. Extensive numerical experiments indicate that Algorithm \ref{1515} performs well, yielding small relative errors compared to M-LBFGS with Armijo line search \cite{Shi2016ALM}. Although Algorithm \ref{1515} theoretically converges to a critical point rather than a global minimum, the comparative results demonstrate its effectiveness in practice, reliably computing the largest eigenvalue with reasonable accuracy. Hence, AggMBFGS, combined with the backtracking modified Armijo line search (Algorithm \ref{1515}), serves as a suitable alternative for computing the largest eigenvalue of high dimensional positive definite matrices.
The worst-case complexity analysis of AggMBFGS for nonconvex problems remains an open problem. 

\subsection*{Competing Interests and Funding}
The author has no competing interests to declare that are relevant to the content of this article. No funding was received to assist with the preparation of this manuscript.
\subsection*{Data Availability}
Data sharing does not apply to this article as no datasets were generated or analyzed during this study.

\bibliography{sn-article}
\end{document}